\documentclass[10pt,twoside]{article}
\usepackage{graphicx}
\usepackage{amsmath}
\usepackage{Latex-document}
\usepackage{amsfonts,amstext,amsbsy,amsopn,amsthm,upref}

\theoremstyle{definition}

\title{\bf  Asymptotic Behaviour of Ergodic\vskip -2mm Integrals of
`Renormalizable' \vskip -2mm Parabolic Flows\vskip 6mm}

\author{G. Forni\vspace*{-0.5cm}\thanks{Department of Mathematics,
Northwestern University, Lunt Hall, 2033 Sheridan Road, Evanston,
IL, 60208-2730, USA. E-mail: gforni@math.northwestern.edu}}
\date{\vspace{-8mm}}

\begin{document}

\maketitle

\thispagestyle{first} \setcounter{page}{317}

\begin{abstract}

\vskip 3mm

Ten years ago A. Zorich discovered, by computer experiments on interval
exchange transformations, some striking new power laws for the ergodic
integrals of generic non-exact Hamiltonian flows on higher genus surfaces.
In Zorich's later work and in a joint paper authored by M. Kontsevich,
Zorich and Kontsevich were able to explain conjecturally most of Zorich's
discoveries by relating them to the ergodic theory of Teichm\"uller flows
on moduli spaces of Abelian differentials.

\smallskip
In this article, we outline a generalization of the Kontsevich-Zorich
framework to a class of `renormalizable' flows on `pseudo-homogeneous'
spaces. We construct for such flows a `renormalization dynamics' on an
appropriate `moduli space', which generalizes the Teichm\"uller flow.
If a flow is renormalizable and the space of smooth functions is `stable',
in the sense that the Lie derivative operator on smooth functions has
closed range, the behaviour of ergodic integrals can be analyzed, at
least in principle, in terms of an Oseledec's decomposition for a
`renormalization cocycle' over the bundle of `basic currents' for the
orbit foliation of the flow.

\smallskip
This approach was suggested by the author's proof of the Kontsevich-Zorich
conjectures and it has since been applied, in collaboration with L. Flaminio,
to prove that the Zorich phenomenon generalizes to several classical examples
of volume preserving, uniquely ergodic, parabolic flows such as horocycle
flows and nilpotent flows on homogeneous $3$-manifolds.

\vskip 4.5mm


\end{abstract}

\vskip 12mm

\section{Introduction} \label{section 1}\setzero

\vskip-5mm

\hspace{5mm}

A fundamental problem in smooth ergodic theory is to establish quantitative
estimates on the asymptotics behaviour of ergodic integrals of smooth
functions. For several examples of {\it hyperbolic }flows, such as geodesic
flows on compact manifolds of negative curvature, the asymptotic behaviour
of ergodic integrals is described by the {\it Central Limit Theorem }(Y.
Sinai, M. Ratner). In these cases, the dynamical system can be described
as an approximation of a `random' stochastic process, like the outcomes of
flipping a coin. Non-hyperbolic systems as not as well understood, with the
important exception of toral flows. For generic non-singular area-preserving
flows on the $2$-torus logarithmic bounds on ergodic integrals of zero-average
functions of bounded variation can be derived by the {\it Denjoy-Koksma
inequality }and the theory of continued fractions. For a general ergodic
flow, ergodic integrals are bounded for all times for a special class of
functions: {\it coboundaries }with bounded `transfer' functions
(Gottschalk-Hedlund). In the hyperbolic examples and in the case of
generic toral flows, a smooth function is a coboundary if and only if it
has zero average with respect to all invariant measures.

\smallskip
In this article, we are interested in flows with {\it parabolic }behaviour.
Following A. Katok, a dynamical system is called parabolic if the rate of
divergence of nearby orbits is at most polynomial in time, while hyperbolic
systems are characterized by exponential divergence. Toral flows are a rather
special parabolic example, called {\it elliptic}, since there is no divergence
of orbits. It has been known for many years that typical examples of parabolic
flows, such as horocycle flows or generic nilpotent flows  are {\it uniquely
ergodic}, but until recently not much was known on the asymptotic behaviour
of ergodic averages, with the exception of some polynomial bounds on the
speed of convergence in the horocycle case (M. Ratner, M. Burger), related
to the polynomial rate of mixing. We have been able to prove, in collaboration
with L. Flaminio, that for many examples of parabolic dynamics the behaviour
of ergodic averages is typically described as follows.

\smallskip
A smooth flow $\Phi^X$ on a finite dimensional manifold $M$ has {\it deviation
spectrum }$\{\lambda_1> ... > \lambda_i> ... >0\}$ with multiplicities $m_1,
...,m_i,...\in {\mathbb Z}^+$ if there exists a system $\{{\mathcal D}_{ij}
\,|\, i\in {\mathbb Z}^+\,,\,\,1\le j\le m_i\}$ of linearly independent
$X$-{\it invariant distributions }such that, for almost all $p\in M$, the
ergodic integrals of any smooth function $f\in C_0^{\infty}(M)$ have an
asymptotic expansion
\begin{equation}
\label{DS}
\int_0^T f\bigl(\Phi^X(t,p)\bigr) dt \,\,= \,\,\sum_{i\in \mathbb{N}}
\sum_{j=1}^{m_i} c_{ij}(p,T) \,{\mathcal D}_{ij}(f)\, T^{\lambda_i}\,\, +\,\,
R(p,T)(f)\,,
\end{equation}
where the real coefficients $c_{ij}(p,T)$ and the distributional remainder
$R(p,T)$ have, for almost all $p\in M$, a sub-polynomial behaviour, in the
sense that
\begin{equation}
\limsup_{T\to +\infty} \frac{ \log \sum_{j=1}^{m_i}|c_{ij}(p,T)|^2}{\log T}\,
=\,\limsup_{T\to +\infty} \frac{\log \|R(p,T)\|}{\log T}\,=\,0\,.
\end{equation}

\smallskip
The notion of a deviation spectrum first arose in the work of A. Zorich
and in his joint work with M. Kontsevich on non-exact Hamiltonian flows with
isolated saddle-like singularities on compact higher genus surfaces. Zorich
discovered in numerical experiments on interval exchange transformations an
unexpected new phenomenon \cite{Zone}. He found that, although a generic flow
on a surface of genus $g\ge 2$ is
uniquely ergodic (H. Masur, W. Veech), for large times the homology classes
of return orbits exhibit unbounded polynomial deviations with exponents
$\lambda_1>\lambda_2>...>\lambda_g>0$ from the line spanned in the homology
group by the {\it Schwartzmann's asymptotic cycle}. In his later work \cite
{Ztwo}, \cite{Zthree} and in joint work with M. Kontsevich \cite{K}, Zorich
was able to explain this phenomenon in terms of conjectures on the Lyapunov
exponents of the Teichm\"uller flow on moduli spaces of holomorphic
differentials on Riemann surfaces. Kontsevich and Zorich also conjectured
that Zorich's phenomenon is not merely topological, but it extends to ergodic
integrals of smooth functions. ``There is, presumably, an equivalent way of
describing the numbers $\lambda_i$. Namely, let $f$ be a smooth function ...
Then for a generic trajectory $p(t)$, we expect that the number $\int_0^T
f\bigl(p(t)\bigr) dt$ for large $T$ with high probability has size
$T^{\lambda_i+o(T)}$ for some $i\in \{1,...,g\}$. The exponent $\lambda_1$
appears for all functions with non-zero average value. The next exponent,
$\lambda_2$, should work for functions in a codimension $1$ subspace of
$C^{\infty}(S)$, etc.'' \cite{K}

\smallskip
Around the same time, we proved that for a generic non-exact Hamiltonian
flow $\Phi^X$ on a higher genus surface not all smooth zero average
functions are smooth coboundaries \cite{Fone}. In fact, we found that, in
contrast with the hyperbolic case and the elliptic case of toral flows,
there are $X$-invariant distributional obstructions, which are not signed
measures, to the existence of smooth solutions of the {\it cohomological
equation }$Xu=f$. This result suggested that Zorich's phenomenon should
be related to the presence of invariant distributions other than the (unique)
invariant probability measure. In fact, in \cite{Ftwo} we were able prove
the Kontsevich-Zorich conjectures that the deviation exponents are non-zero
and that generic non-exact Hamiltonian flows on higher genus surfaces have
a deviation spectrum. Recently, in collaboration with L. Flaminio, we have
proved that other classical parabolic examples, such as horocycle flows on
compact surfaces of constant negative curvature \cite{FFone} and generic
nilpotent flows on compact $3$-dimensional nilmanifolds \cite{FFtwo}, do
have a deviation spectrum, but of countable multiplicity, in contrast with
the case of flows on surfaces which have spectrum of finite multiplicity
equal to the genus.

\smallskip
We will outline below a general framework, derived mostly from \cite{Ftwo}
and successfully carried out in \cite{FFone}, \cite{FFtwo}, for proving that
a flow on a {\it pseudo-homogeneous space }has a deviation spectrum. Our
framework is based on the construction of an appropriate {\it renormalization
dynamics }on a moduli space of pseudo-homogeneous structures, which generalizes
the Teichm\"uller flow. A {\it renormalizable flow }for which the space of
smooth functions is {\it stable} (in the sense of A. Katok), has a deviation
spectrum determined by the Lyapunov exponents of a {\it renormalization
cocycle }over a bundle of {\it basic currents}. Pseudo-homogeneous spaces
are a generalization of homogeneous spaces. The motivating non-homogeneous
example is given by any punctured Riemann surface carrying a holomorphic
differential vanishing only at the punctures. It turns out that renormalizable
flows are necessarily parabolic. In fact, the class of renormalizable flows
encompasses all parabolic flows which are reasonably well-understood, while
not much is known for most non-renormalizable parabolic flows, such as generic
geodesic flows on flat surfaces with conical singularities.  Our approach
unifies and generalizes several classical quantitative equidistribution
results such as the Zagier-Sarnak results for periodic horocycles on
non-compact hyperbolic surfaces of finite volume \cite{FFone} or number
theoretical results on the asymptotic behaviour of theta sums \cite{FFtwo}.

\section{Renormalizable flows} \label{section 2}
\setzero\vskip-5mm \hspace{5mm }

Let $\mathfrak g$ be a finite dimensional real Lie algebra. A $\mathfrak
g$-{\it structure }on a manifold $M$ is defined to be a homomorphism
$\tau$ from $\mathfrak g$ into the Lie algebra ${\mathcal V}(M)$ of all
smooth vector fields on $M$. This notion is well-known in the theory of
transformation groups (originated in the work of S. Lie) under the name of
`infinitesimal $G$-transformation group' (for a Lie group $G$ with $\mathfrak
g$ as Lie algebra). The second fundamental theorem of Lie states that any
infinitesimal $G$-transformation group $\tau$ on $M$ can be `integrated'
to yield an essentially unique local $G$-transformation group.
A $\mathfrak g$-structure $\tau$ will be called {\it faithful }if $\tau$
induces a linear isomorphism from ${\mathfrak g}$ onto $T_x M$, for all $x\in
M$. Let $\tau$ be a $\mathfrak g$-structure. For each element $X\in {\mathfrak
g}$, the vector field $X_{\tau}:=\tau(X)$ generates a (partially defined) flow
$\Phi^X_{\tau}$ on $M$. Let $E_t(X_{\tau})\subset M$ be the closure of the
complement of the domain of definition of the map $\Phi^X_{\tau}(t,\cdot)$
at time $t\in {\mathbb R}$. A faithful $\mathfrak g$-structure will be called
{\it pseudo-homogeneous }if for every $X\in \mathfrak g$ there exists $t>0$
such that $E_t(X_{\tau}) \cup E_{-t}(X_{\tau})$ has zero (Lebesgue) measure.
A manifold $M$ endowed with a pseudo-homogeneous $\mathfrak g$-structure will
be called a {\it pseudo-homogeneous }$\mathfrak g$-{\it space}. All
homogeneous spaces are pseudo-homogeneous.

\smallskip
Let ${\mathcal T}_{\mathfrak g}(M)$ be the space of all pseudo-homogeneous
$\mathfrak g$-structures on $M$. The automorphism group $\text{Aut}(\mathfrak
g)$ acts on ${\mathcal T}_{\mathfrak g}(M)$ by composition on the right. The
group $\text{Diff}(M)$ acts on ${\mathcal T}_{\mathfrak g}(M)$ by composition
on the left. The spaces
\begin{equation}
\text{T}_{\mathfrak g}(M):={\mathcal T}_{\mathfrak g}(M)/\text{Diff}_0(M)
\,\,,\,\,\,\,
{\mathcal M}_{\mathfrak g}(M):=\text{T}_{\mathfrak g}(M)/\Gamma(M) \,,
\end{equation}
where $\Gamma(M):=\text{Diff}^+(M)/\text{Diff}_0(M)$ is the {\it mapping
class group}, will be called respectively the {\it Teichm\"uller space }and
the {\it moduli space }of pseudo-homogeneous $\mathfrak g$-structures on $M$.
The group $\text{Aut}(\mathfrak g)$ acts on the Teichm\"uller space $\text{T}
_{\mathfrak g}(M)$ and on the moduli space ${\mathcal M}_{\mathfrak g}(M)$,
since in both cases the action of $\text{Aut}(\mathfrak g)$ on  ${\mathcal
T}_{\mathfrak g}(M)$ passes to the quotient.

\smallskip
Let $\text{Aut}^{(1)}(\mathfrak g)$ be the subgroup of automorphisms with
determinant one. An element $X\in \mathfrak g$ will be called {\it a priori
renormalizable }if there exists a partially hyperbolic one-parameter subgroup
$\{G^X_t\}\subset \text{Aut}^{(1)}(\mathfrak g)$, $t\in {\mathbb R}$ ($t\in
{\mathbb Z}$), in general non-unique, with a single (simple) Lyapunov exponent
$\mu_X>0$ such that
\begin{equation}
\label{eq:RN}
G^X_t(X)=e^{t\mu_X}\,X \,\,.
\end{equation}
It follows from the definition that the subset of a priori renormalizable
elements of a Lie algebra $\mathfrak g$ is saturated with respect to the
action of $\text{Aut}(\mathfrak g)$. The subgroup $\{G^X_t\}$ acts on the
Teichm\"uller space and on the moduli space of pseudo-homogeneous $\mathfrak
g$-structures as a `renormalization dynamics' for the family of flows
$\Phi^X_{\tau}$ generated by the vector fields $\{X_{\tau}\,|\,\tau\in
{\mathcal T}_{\mathfrak g}(M)\}$ on $M$. It will be called a {\it generalized
Teichm\"uller flow (map)}. A flow $\Phi^X_{\tau}$ will be called {\it
renormalizable }if $\tau\in {\mathcal M}_{\mathfrak g}(M)$ is a recurrent
point for some generalized Teichm\"uller flow (map) $G_t^X$. If $\mu$ is
a probability $G^X_t$-invariant measure on the moduli space, then by
Poincar\'e recurrence the flow $\Phi^X_{\tau}$ is renormalizable for
$\mu$-almost all $\tau\in {\mathcal M}_{\mathfrak g}(M)$.

\smallskip
Let $R$ be an inner product on $\mathfrak g$. Every faithful $\mathfrak g$
structure $\tau$ induces a Riemannian metric $R_{\tau}$ of constant curvature
on $M$. Let $\omega_{\tau}$ be the volume form of $R_{\tau}$. The total volume
function $A:{\mathcal T}_{\mathfrak g}(M)\to {\mathbb R}^+ \cup \{+\infty\}$
is $\text{Diff}^+(M)$-invariant and $\text{Aut}^{(1)}(\mathfrak g)$-invariant.
Hence $A$ is well-defined as an $\text{Aut}^{(1)}(\mathfrak g)$-invariant
function on the Teichm\"uller space and on the moduli space. It follows that
the subspace of finite-volume $\mathfrak g$-structures has an $\text{Aut}^{(1)}
(\mathfrak g)$-invariant stratification by the level hypersurfaces of the total
volume function. Since different hypersurfaces are isomorphic up to a dilation,
when studying finite-volume spaces it is sufficient to consider the
hypersurface of volume-one $\mathfrak g$-structures:
\begin{equation}
\text{T}_{\mathfrak g}^{(1)}(M):=\text{T}_{\mathfrak g}(M)\cap A^{-1}(1)
\,\,,\,\,\,\,
{\mathcal M}_{\mathfrak g}^{(1)}(M):= {\mathcal M}_{\mathfrak g}(M)
\cap A^{-1}(1) \,.
\end{equation}
Let $\tau$ be a faithful $\mathfrak g$-structure and let $X\in\mathfrak g$. If
the linear map $ad_X$ on $\mathfrak g$ has zero trace, the flow $\Phi^X_{\tau}$
preserves the volume form $\omega_{\tau}$ and $X_{\tau}$ defines a symmetric
operator on $L^2(M,\omega_{\tau})$ with domain $C_0^{\infty}(M)$. If $\tau$
is pseudo-homogeneous, by E. Nelson's criterion \cite{N}, $X_{\tau}$ is
essentially skew-adjoint. It turns out that any a priori renormalizable
element $X\in \mathfrak g$ is {\it nilpotent}, in the sense that all
eigenvalues of the linear map $ad_X$ are equal to zero, hence the flow
$\Phi^X_{\tau}$ is volume preserving and parabolic. In all the examples
we have considered, the Lie algebra $\mathfrak g$ is {\it traceless}, in
the sense that for every element $X\in \mathfrak g$, the linear map $\text
{ad}_X$ has vanishing trace. In this case, any pseudo-homogeneous $\mathfrak
g$-structure induces a representation of the Lie algebra $\mathfrak g$ by
essentially skew-adjoint operators on the Hilbert space $L^2(M,\omega_{\tau})$
with common invariant domain $C_0^{\infty}(M)$.

\section{Examples} \label{section 3}
\setzero\vskip-5mm \hspace{5mm }

Homogeneous spaces provide a wide class of examples. Let $G$ be a finite
dimensional (non-compact) Lie group with Lie algebra $\mathfrak g$ and
let $M=G/\Gamma$ be a (compact) homogeneous space. The Teichm\"uller
space $\text{T}_G(M)\subset \text{T}_{\mathfrak g}(M)$ and the moduli space
${\mathcal M}_G(M)\subset {\mathcal M}_{\mathfrak g}(M)$ of all homogeneous
$G$-space structures on $M$ are respectively isomorphic to the Lie group
$\text{Aut}(G)$ and to the homogeneous space $\text{Aut}(G)/\text{Aut}(G,
\Gamma)$, where $\text{Aut}(G,\Gamma)<\text{Aut}(G)$ is the subgroup of
automorphisms which stabilize the lattice $\Gamma$. The Teichm\"uller
and moduli spaces $\text{T}^{(1)}_G(M)$ and ${\mathcal M}^{(1)}_G(M)$ of
homogeneous volume-one $G$-space structures are respectively isomorphic
to the subgroup $\text{Aut}^{(1)}(G)$ of orientation preserving, volume
preserving automorphisms and to the homogeneous space $\text{Aut}^{(1)}(G)
/\text{Aut}^{(1)}(G,\Gamma)$.

\smallskip
In the {\it Abelian }case ${\mathfrak g}={\mathbb R}^n$, any $X\in {\mathbb
R}^n$ is a priori renormalizable. In fact, the group $\text{Aut}^{(1)}
({\mathbb R}^n)=\text{SL}(n,{\mathbb R})$ acts transitively on ${\mathbb R}^n$
and $X_1=(1,0,...,0)$ is renormalized by the one-parameter group $G_t:=\text
{diag}(e^t,e^{-t/n},...,e^{-t/n})\subset \text{SL}(n,{\mathbb R})$.
Finite volume Abelian homogeneous spaces are diffeomorphic to $n$-dimensional
tori ${\mathbb T}^n$. The generalized Teichm\"uller flow $G_t$ on the moduli
space of all volume-one Abelian homogeneous structures on ${\mathbb T}^n$ is
a volume preserving Anosov flow on the finite-volume non-compact manifold
$\text{SL}(n,\mathbb R)/\text{SL}(n,\mathbb Z)$. Hence, in this case, almost
all homogeneous flows are renormalizable, by Poincar\'e recurrence theorem.
The dynamics of the flow $G_t$ has been investigated in depth by D. Kleinbock
and G. Margulis in connection with the theory of Diophantine approximations.

\smallskip
In the {\it semi-simple }case, let ${\mathfrak g}={\mathfrak s}{\mathfrak l}
(2,\mathbb R)$ be the unique  $3$-dimensional simple Lie algebra. There is a
basis $\{H,H^{\perp},X\}$ with commutation relations $[X,H]=H$, $[X,H^{\perp}]
=-H^{\perp}$ and $[H,H^{\perp}]=2X$. The elements $H$, $H^{\perp}$ are
renormalized by the one-parameter group $G_t:=\text{diag}(e^t,e^{-t},1)\subset
\text{Aut}^{(1)}(\mathfrak g)$, while $X$ is not a priori renormalizable. The
unit tangent bundle of any hyperbolic surface $S$ can be identified to a
homogeneous $\mathfrak g$-space $M:=\text {PSL}(2,{\mathbb R})/\Gamma$. The
vector fields $H$, $H^{\perp}\in {\mathfrak g}$ generate the horocycle flows
and the vector field $X$ generates the geodesic flow on $S$. Since $G_t$ is a
group of {\it inner }automorphisms, it is in fact generated by the geodesic
vector field $X$, every point of the moduli space is fixed under $G_t$. Hence
horocycle flows are renormalizable on every homogeneous space $\text{PSL}
(2,{\mathbb R})/\Gamma$.

\smallskip
In the {\it nilpotent}, non-Abelian case, let $\mathfrak n$ be the Heisenberg
Lie algebra, spanned by elements $\{X,X^{\perp},Z\}$ such that $[X,X^{\perp}]
=Z$ and $Z$ is a generator of the one-dimensional center $Z_{\mathfrak n}$.
The element $X$ is renormalized by the one-parameter subgroup $G_t:=\text{diag}
(e^t,e^{-t},1)\subset \text{Aut}^{(1)}({\mathfrak n})$. Since the group
$\text{Aut}({\mathfrak n})$ acts transitively on ${\mathfrak n}\setminus
Z_{\mathfrak n}$, every $Y\in{\mathfrak n}\setminus Z_{\mathfrak n}$ is a
priori renormalizable, while the elements of the center are not. A compact
nilmanifold modeled over the Heisenberg group $N$ is a homogeneous space
$M=N/\Gamma$, where $\Gamma$ is a co-compact lattice. These spaces are
topologically circle bundles over ${\mathbb T}^2$ classified by their Euler
characteristic. The moduli space ${\mathcal M}^{(1)}_N(M)$ of volume-one
homogeneous $\mathfrak n$-structures on $M$ is a $5$-dimensional finite-volume
non-compact orbifold which fibers over the modular surface $\text{SL}(2,
{\mathbb R})/\text{SL}(2,{\mathbb Z})$ with fiber ${\mathbb T}^2$. The
generalized Teichm\"uller flow is an Anosov flow on ${\mathcal M}^{(1)}_N(M)$
\cite{FFtwo}.

\smallskip
The motivation for our definition of a pseudo-homogeneous space comes
from the theory of Riemann surfaces of higher genus. Any holomorphic
(Abelian) differential $h$ on a Riemann surface $S$ of genus $g\ge 2$,
vanishing at $Z_h\subset S$, induces a (non-unique) pseudo-homogeneous
${\mathbb R}^2$-structure on the open manifold $M_h:=S\setminus Z_h$
In fact, the frame $\{X,X^{\perp}\}$ of $TS|M_h$ uniquely determined by
the conditions
\begin{equation}
\frac{\sqrt{-1}}{2}\, \imath_X\,(h\wedge {\bar h})\,=\, \Im(h) \,\,\,\,,
\,\,\,\,
\frac{\sqrt{-1}}{2}\, \imath_{X^{\perp}} \,(h\wedge {\bar h})\,=\,-\Re(h)
\end{equation}
satisfies the Abelian commutation relation $[X,X^{\perp}]=0$ and the
homomorphism $\tau_h:{\mathbb R}^2\to {\mathcal V}(M_h)$ such that $\tau_h
(1,0)=X$, $\tau_h(0,1)=X^{\perp}$ is a pseudo-homogeneous ${\mathbb
R}^2$-structure on $M_h$. Let $Z\subset S$ be a given subset of cardinality
$\sigma\in{\mathbb N}$ and let $\kappa=(k_1,...,k_{\sigma})\in ({\mathbb
Z}^+)^{\sigma}$ with $\sum k_i=2g-2$. Let ${\mathcal H}_{\kappa}(S,Z)$ be
the space of Abelian differentials $h$ with $Z_h=Z$ and zeroes of
multiplicities $(k_1,...,k_{\sigma})$. The projection of the set
$\{\tau_h\in {\mathcal T}_{{\mathbb R}^2}(M)\,|\,h\in {\mathcal H}_{\kappa}
(S,Z)\}$ into the moduli space ${\mathcal M}_{{\mathbb R}^2}(M)$ of
pseudo-homogeneous ${\mathbb R}^2$-structures on $M:=S\setminus Z$ is
isomorphic to a stratum ${\mathcal H}({\kappa})$ of the moduli space of
Abelian differentials on $S$. The flow induced on ${\mathcal H}({\kappa})$
by the one-parameter group of automorphism $G_t=\text{diag}(e^t,e^{-t})
\subset \text{SL}(2,\mathbb R)$ coincides with the Teichm\"uller flow on
the stratum ${\mathcal H}({\kappa})$.

\section{Cohomological equations} \label{section 4}
\setzero\vskip-5mm \hspace{5mm }

Let $(\mathfrak g, R)$ be a finite dimensional Lie algebra endowed with an
inner product. Any pseudo-homogeneous $\mathfrak g$-structure $\tau$ on
a manifold $M$ induces a {\it Sobolev filtration }$\{\text{W}^s_{\tau}
(M)\}_{s\ge 0}$ on the space $W^0_{\tau}(M):=L^2(M,\omega_{\tau})$ of
square-integrable functions. Let $\triangle_{\tau}$ be the non-negative
{\it Laplace-Beltrami
operator }of the Riemannian metric $R_{\tau}$ on $M$. The Laplacian is densely
defined and symmetric on the Hilbert space $W^0_{\tau}(M)$ with domain
$C_0^{\infty}(M)$, but it is not in general essentially self-adjoint. In fact,
if $\mathfrak g$ is traceless, by a theorem of E. Nelson \cite{N}, $\triangle
_{\tau}$ is essentially self-adjoint if and only if the representation $\tau$
of the Lie algebra $\mathfrak g$ on $W^0_{\tau}(M)$ by essentially skew-adjoint
operators induces a unitary representation of a Lie group. Let then $\bar
\triangle_{\tau}$ be the {\it Friederichs extension }of $\triangle_{\tau}$.
The Sobolev space $W^s_{\tau}(M)$, $s>0$, is defined as the maximal domain
of the operator $(\text{I}+{\bar\triangle}_{\tau})^{s/2}$ endowed with the
norm

\begin{equation}
\|f\|_{s,\tau}:= \|(\text{I}+\bar{\triangle})^{s/2}f\|_{0,\tau}\,.
\end{equation}

The Sobolev spaces $W^{-s}_{\tau}(M)$ are defined as the duals of the Hilbert
spaces $W^s_{\tau}(M)$, for all $s > 0$. Let $C_B^0(M)$ be the space of
continous bounded functions on $M$. The pseudo-homogeneous space $(M,\tau)$
will be called of {\it bounded type }if there is a continous (Sobolev) embedding
$W^s_{\tau}(M)\subset C_B^0(M)$ for all $s>\text{dim}(M)/2$. The bounded-type
condition is essentially a geometric property of the pseudo-homogeneous
structure.

\smallskip
Let $X\in {\mathfrak g}$. Following A. Katok, the space $W^s_{\tau}(M)$ is
called $W^t_{\tau}(M)$-{\it stable }with respect to the flow $\Phi^X_{\tau}$
if the subspace
\begin{equation}
R^{s,t}(X_{\tau}):= \{ f\in W^s_{\tau}(M)\,|\, f=X_{\tau}u\,\,,\,\,\,\,
u\in W^t_{\tau}(M)\}
\end{equation}
is closed in $W^s_{\tau}(M)$. The flow $\Phi^X_{\tau}$ will be called {\it
tame }(of degree $\ell>0$) if $W^s_{\tau}(M)$ is $W^{s-\ell}_{\tau}(M)$-stable
with respect to $\Phi^X_{\tau}$ for all $s>\ell$. In all the examples of \S 3,
generic renormalizable flows are tame. In particular, it is well known that
generic toral flows are tame, horocycle flows and generic nilpotent flows on
$3$-dimensional compact nilmanifolds were proved tame of any degree $\ell>1$
in \cite{FFone}, \cite{FFtwo}, generic non-exact Hamiltonian flows on higher
genus surfaces were proved tame in \cite{Fone}. These results are based on the
appropriate harmonic analysis: in the homogeneous cases, the theory of unitary
representations for the Lie group $\text{SL}(2,{\mathbb R})$ \cite{FFone} and
the Heisenberg group \cite{FFtwo}; in the more difficult non-homogeneous case
of higher genus surfaces, the theory of boundary behaviour of holomorphic
functions on the unit disk plays a crucial role \cite{Ftwo}.

\smallskip
If the Sobolev space $W^s_{\tau}(M)$ is stable with respect to the flow
$\Phi^X_{\tau}$, the closed range $R^{s,t}(X_{\tau})$ of the operator
$X_{\tau}$ coincides with the distributional kernel ${\mathcal I}^s(X_{\tau})
\subset W^{-s}_{\tau}(M)$ of $X_{\tau}$, which is a space of $X_{\tau}$-{\it
invariant distributions}. Let $X$ be any smooth vector field on a manifold
$M$. A distribution ${\mathcal D}\in{\mathcal D}'(M)$ is called $X$-invariant
if $X{\mathcal D}=0$ in ${\mathcal D}'(M)$. Invariant distributions are in
bijective correspondence with (homogeneous) one-dimensional {\it basic
currents }for the orbit foliation ${\mathcal F}(X)$ of the flow $\Phi^X$. A
one-dimensional basic current $C$ for a foliation $\mathcal F$ on $M$ is a
continous linear functional on the space $\Omega^1_0(M)$ of smooth $1$-forms
with compact support such that, for all vector fields $Y$ tangent to
$\mathcal F$,
\begin{equation}
\label{eq:BC}
\imath_YC\,=\,{\mathcal L}_YC\,=\,0\,\, (\iff \,\, \imath_YC\,=\,dC\,=\,0)\,.
\end{equation}
It follows from the definitions that the one-dimensional current $C:=\imath_X
{\mathcal D}$ is basic for ${\mathcal F}(X)$ if and only if the distribution
$\mathcal D$ is $X$-invariant. Let ${\mathcal I}(X)$ be the space of all
$X$-invariant distributions and ${\mathcal B}(X)$ be the space of all
one-dimensional basic currents for the orbit foliation ${\mathcal F}(X)$. The
linear map $\imath_X:{\mathcal I}(X)\to {\mathcal B}(X)$ is bijective.

\smallskip
Let $(M,\tau)$ be a pseudo-homogeneous space. There is a
well-defined Hodge (star) operator and a space $\text{C}^0_{\tau}(M)$ of
square-integrable $1$-forms on $M$ associated with the metric $R_{\tau}$.
Since the Laplace operator $\triangle_{\tau}$ extends to $\text{C}^0_{\tau}
(M)$ with domain $\Omega^1_0(M)$, it is possible to define, as in the case
of functions, a Sobolev filtration $\{\text{C}^s_{\tau}(M)\}_{s\ge 0}$, on
the space $\text{C}^0_{\tau}(M)$. The Sobolev spaces $\text{C}^{-s}_{\tau}
(M)$ are defined as the duals of the Sobolev spaces $\text{C}^s(M)$, for
all $s>0$. Let ${\mathcal B}^s(X_{\tau}):={\mathcal B}(X_{\tau})\cap
\text{C}^{-s}_{\tau}(M)$ be the subspaces of basic currents of Sobolev
order $\le s$ for the orbit foliation ${\mathcal F}(X_{\tau})$. The
space ${\mathcal B}^s(X_{\tau})$ is the image of ${\mathcal I}^s(X_{\tau})
:={\mathcal I}(X_{\tau})\cap W^{-s}_{\tau}(M)$ under the bijective map
$\imath_X:{\mathcal I}(X)\to {\mathcal B}(X)$. In the case of minimal
toral flows the space ${\mathcal B}^s(X_{\tau})$ is one-dimensional for
all $s\ge 0$ (as all invariant distributions are scalar multiples of the
unique invariant probability measure). In the parabolic examples we
have studied, ${\mathcal B}^s(X_{\tau})$ has countable dimension,
as soon as $s>1/2$, for horocycle flows or generic nilpotent flows, while for
generic non-exact Hamiltonian flows on higher genus surfaces the dimension is
finite for all $s>0$ and grows linearly with respect to $s>0$. This finiteness
property seems to be an exceptional low dimensional feature.

\section{The renormalization cocycle} \label{section 5}
\setzero\vskip-5mm \hspace{5mm }

The Sobolev spaces $\text{C}^s_{\tau}(M)$ of one-dimensional currents form a
smooth infinite dimensional vector bundle over ${\mathcal T}_{\mathfrak g}(M)$.
Such bundles can be endowed with a flat connection with parallel transport
given locally by the identity maps $\text{C}^s_{\tau}(M)\to\text{C}^s_{\tau'}
(M)$, for any $\tau\approx\tau'\in {\mathcal T}_{\mathfrak g}(M)$. Since the
diffeomorphism group $\text{Diff}(M)$ acts on $\text{C}^s_{\tau}(M)$ by
push-forward, we can define (orbifold) vector bundles $\text{C}^s_{\mathfrak g}
(M)$ over the Teichm\"uller space $\text{T}_{\mathfrak g}(M)$ or the moduli
space ${\mathcal M}_{\mathfrak g}(M)$ of pseudo-homogeneous structures on $M$.
If $X\in \mathfrak g$ is a priori renormalizable, a generalized Teichm\"uller
flow (map) $G^X_t$ can be lifted by parallel transport to a `renormalization
cocycle' $R^X_t$ on the bundles of currents $\text{C}^s_{\mathfrak g}(M)$ over
the Teichm\"uller space or the moduli space. It follows from the definitions
that the sub-bundles ${\mathcal B}^s_{\mathfrak g}(X)\subset \text{C}^s_
{\mathfrak g}(M)$ with fibers the subspaces of basic currents ${\mathcal B}^s
(X_{\tau})\subset \text{C}^s_{\tau}(M)$ are $R^X_t$-invariant. It can be proved
that, for any $G_t$-ergodic probability measure $\mu$ on the moduli space, if
the flows $\Phi^X_{\tau}$ are tame of degree $\ell>0$ for $\mu$-almost all
$\tau\in{\mathcal M}_{\mathfrak g}(M)$, then the sub-bundles ${\mathcal B}^s
_{\mathfrak g}(X)$ are $\mu$-almost everywhere defined with closed (Hilbert)
fibers of constant rank, for all $s>\ell$.

\smallskip
In the examples considered, with the exception of flows on higher genus
surfaces, the Hilbert bundles of basic currents ${\mathcal B}^s_{\mathfrak g}
(X)$ are infinite dimensional, and to the author's best knowledge, available
Oseledec-type theorems for Hilbert bundles do not apply to the renormalization
cocycle. However, the cocycle has a well defined Lyapunov spectrum and an
Oseledec decomposition. We are therefore led to formulate the following
hypothesis:

\smallskip
\noindent $H_1(s)$.  The renormalization cocyle $R^X_t$ on the bundle
${\mathcal B}^s_{\mathfrak g}(X)$ over the dynamical system $(G_t^X,\mu)$
has a Lyapunov spectrum $\{\nu_1> ...> \nu_k> ...>0>...\}$ and an Oseledec's
decomposition
\begin{equation}
\label{eq:OD}
{\mathcal B}^s_{\mathfrak g}(X)=E^s_{\mathfrak g}(\nu_1)\oplus ... \oplus
E^s_{\mathfrak g}(\nu_k)\oplus ... \oplus N^s_{\mathfrak g}\,\,,
\end{equation}
in which the components $E^s_{\mathfrak g}(\nu_k)$ correspond to the Lyapunov
exponents $\nu_k>0$, while the component $N^s_{\mathfrak g}$ has a non-positive
top Lyapunov exponent. Our result on the existence of a deviation spectrum
requires an additional technical hypothesis, verified in our examples.

\smallskip
\noindent $H_2(s)$. Let $\gamma_{\tau}^1(p)$ be the one-dimensional current
defined by the time $T=1$ orbit-segment of the flow $\Phi^X_{\tau}$ with
initial point $p\in M$. $(a)$ The essential supremum of the norm $\|\gamma
_{\tau}^1(p)\|_{\tau,s}$ over $p\in M$ is locally bounded for $\tau\in\text
{supp}(\mu)\subset {\mathcal M}_{\mathfrak g}(M)$; $(b)$ The orthogonal
projections of $\gamma_{\tau}^1(p)$ on all subspaces $E^s_{\tau}(\nu_k)
\subset \text{C}^{-s}_{\mathfrak g}(M)$ are non-zero for $\mu$-almost all
$\tau \in {\mathcal M}_{\mathfrak g}(M)$ and almost all $p\in M$.

\smallskip
{ Let }$X\in {\mathfrak g}$ { be a priori renormalizable and let
}$\mu$ { be a }$G_t^X$-{ invariant Borel probability measure on
}${\mathcal M} _{\mathfrak g}(M)$, { supported on a stratum of
bounded-type }$\mathfrak g$-{ structures. If the flow
}$\Phi_{\tau}^X$ { is tame of degree }$\ell
>0$ { and the hypoteses }$H_1(s)$, $H_2(s)$ { are verified for }$s>\ell+
\text{dim}(M)/2$, { for }$\mu$-{ almost }$\tau\in {\mathcal
M}_{\mathfrak g}$, { the flow }$\Phi_{\tau}^X$ { has a deviation
spectrum with deviation exponents}
\begin{equation}
\nu_1/\mu_X \,>\,...\,>\,\nu_k/\mu_X \,>\,...\,>\,0
\end{equation}
{ and multiplicities given by the decomposition }(\ref{eq:OD}) {
of the renormalization coycle.}

\smallskip
In the homogeneous examples, the Lyapunov spectrum of the renormalization
cocycle is computed explicitly in every irreducible unitary representation
of the structural Lie group. In the horocycle case, the existence of an
Oseledec's decomposition (\ref{eq:OD}) is equivalent to the statement that
the space of horocycle-invariant distributions is spanned by (generalized)
eigenvectors of the geodesic flow, well-known in the representation theory
of semi-simple Lie groups as {\it conical distributions }\cite{H}. In the
non-homogeneous case of higher genus surfaces, the Oseledec's theorem applies
since the bundles ${\mathcal B}^s_{\mathfrak g}(X)$ are finite dimensional.
We have found in all examples a surprising heuristic relation between the
Lyapunov exponents of the renormalization cocycle and the Sobolev regularity
of basic currents (or equivalently of invariant distributions): the subspaces
$E^s_{\tau}(\nu_k)$ are generated by basic currents of {\it Sobolev order }$1
-\nu_k/\mu_X\ge 0$. The Sobolev order of a one-dimensional current $C$ is
defined as the infimum of all $s>0$ such that $C\in \text{C}^{-s}_{\tau}(M)$.

\smallskip
In the special case of non-exact Hamiltonian flow on higher genus surfaces
the Lyapunov exponents of the renormalization cocycle are related to those of
the Teichm\"uller flow. In fact, let $S$ be compact orientable surface of
genus $g\ge 2$ and let ${\mathcal H}({\kappa})$ be a stratum of Abelian
differentials vanishing at $Z\subset S$. Let ${\mathcal B}_{\kappa}(X)
\subset {\mathcal B}_{{\mathbb R}^2}(X)$ be the measurable bundle of basic
currents over ${\mathcal H}({\kappa})\subset {\mathcal M}_{{\mathbb R}^2}
(S\setminus Z)$ and let $H^1_{\kappa}(S\setminus Z,\mathbb R)$ be the bundle
over ${\mathcal H}({\kappa})$ with fibers isomorphic to the real cohomology
$H^1(S\setminus Z,\mathbb R)$. Since basic currents are closed, there exists
a {\it cohomology map }$j_{\kappa}:{\mathcal B}_{\kappa}(X) \to H^1_{\kappa}
(S\setminus Z,{\mathbb R})$ such that, as proved in \cite{Ftwo}, the
restrictions $j_{\kappa}|{\cal B}^s_{\kappa}(X)$ are surjective for all
$s>>1$ and, for all $s\ge 1$, there are exact sequences
\begin{equation}
0\rightarrow \mathbb R \rightarrow {\mathcal B}_{\kappa}^{s-1}(X)\,\,
^{\underrightarrow{\,\,\,\, \delta_{\kappa}\,\,\,\,}}\,\, {\mathcal B}
_{\kappa}^s(X)\,\, ^{\underrightarrow{\,\,\,\, j_{\kappa} \,\,\,\,}}\,\,
H^1_{\kappa}(S\setminus Z,{\mathbb R})\,\,.
\end{equation}
The renormalization cocycle $R^X_t$ on ${\mathcal B}_{\kappa}^s(X)$ projects
for all $s>>1$ onto a cocycle on the cohomology bundle $H^1_{\kappa}(S\setminus
Z,{\mathbb R})$, introduced by M. Kontsevich and A. Zorich in order to explain
the {\it homological }asymptotic behaviour of orbits of the flow $\Phi_{\tau}^X$
for a generic $\tau\in {\cal H}(\kappa)\subset {\cal M}_{{\mathbb R}^2}(S
\setminus Z)$ \cite{K}. The Lyapunov exponents of the Kontsevich-Zorich cocycle
on $H^1_{\kappa}(S\setminus Z,{\mathbb R})$,
\begin{equation}
\label{eq:RKZE}
\lambda_1=1>\lambda_2\ge \dots \ge\lambda_g\ge \overbrace{0=\cdots=
0}^{\#Z-1} \ge -\lambda_g\ge \dots \ge -\lambda_2>-\lambda_1=-1 \,,
\end{equation}
are related the Lyapunov exponents of the Teich\"uller flow on
${\mathcal H}({\kappa})$ \cite{K}, \cite{Ftwo}. Since the bundle map
$\delta_{\kappa}$ shifts Lyapunov exponents by $-1$ and, as conjectured
in \cite{K} and proved in \cite{Ftwo}, the Kontsevich-Zorich exponents
$\lambda_1=1>\lambda_2\ge \dots \ge\lambda_g$ are non-zero, the strictly
positive exponents of the renormalization cocycle coincide with the
Kontsevich-Zorich exponents. This reduction explains why in the case of
non-exact Hamiltonian flows on surfaces the Lyapunov exponents of the
Teichm\"uller flow are related to the deviation exponents for the ergodic
averages of smooth functions.

\label{lastpage}

\end{document}